\documentclass[11pt]{article}
\usepackage{amsmath, amssymb}

\usepackage{color}

\newcommand{\R}{\mathbb{R}}

\newcommand{\N}{\mathbb{N}}
\newcommand{\bea}{\begin{eqnarray}}
\newcommand{\eea}{\end{eqnarray}}

\def\de{\delta}
\def\e{\varepsilon}

\def\s{\sigma}

\def\supp{{\rm supp}}

\def\1{\rm Id}
\def\sgn{\rm sgn}
\def\sup{{\rm sup}}

\def\ci{\circ}
\def\cl{{\rm cl}}

\newcommand{\qed}{$\hfill\blacksquare$}

\def\V{\noindent}

\def\Int{{\rm int}}
\def\cl{{\rm cl}}

\newcommand{\bean}{\begin{eqnarray*}}
\newcommand{\eean}{\end{eqnarray*}}

\newtheorem{Theorem}{Theorem}

\newcommand{\ben}{\begin{enumerate}}
\newcommand{\een}{\end{enumerate}}
\newcommand{\bit}{\begin{itemize}}
\newcommand{\eit}{\end{itemize}}
\newcommand{\edoc}{\end{document}}

\usepackage{hyperref}

\usepackage{enumerate}

\title{Maximality and Cauchy developments for Lorentzian length spaces}

\begin{document}

\author{Olaf M\"uller\footnote{Institut f\"ur Mathematik, Humboldt-Universit\"at zu Berlin, Unter den Linden 6, D-10099 Berlin, \texttt{Email: o.mueller@hu-berlin.de}}}

\date{\today}
\maketitle

\begin{abstract}
\V This article suggests the definition of "Lorentzian space" weakening the notion of Lorentzian length spaces just as much that it allows for a functor from the category of strongly causal Lorentzian manifolds to the corresponding category of Lorentzian spaces, and considers three problems in the context of maximal Cauchy developments of Lorentzian spaces: The first is to define pointed Gromov-Hausdorff metrics for spatially and temporally noncompact Lorentzian spaces, the second to present an explicit non-spacetime example of a maximal globally hyperbolic Lorentzian space, the third to define canonical representatives for Cauchy developments. \\ A certain well-posedness property for geodesics plays a key role in each of the problems.   
\end{abstract}

For general relativistic theories like Einstein-Maxwell theory, the notion of 'global solution' is intricate, due to diffeomorphism-invariance. Their domain of definition varies: how to compare solutions in size? Each local solution can be "stretched" to be defined for all times. 
The seminal article to solve these conceptual problems and to define the maximal Cauchy development (mCd) was by Choquet-Bruhat and Geroch \cite{CBG}. Its two main ingredients were harmonic coordinates (in which the equations take a quasilinear symmetric hyperbolic form, thus satisfying a local uniqueness statement allowing for pasting together local solutions), and Zorn's lemma for the partial order of the solutions in this gauge. There are many alternative versions and adaptations of this procedure, notably two approaches \cite{wW}, \cite{jS} of "dezornification" using only weaker versions of Zorn's lemma.

This article makes first steps towards defining mCds for Lorentzian length spaces (LLSs). In LLSs, vastly generalizing $C^1$ causal spacetimes, we have a notion of global hyperbolicity \cite{KS} implying existence of Cauchy sets \cite{BGH} and of Cauchy time functions adapted to given such \cite{oM-GH}, and functorial pseudo-metrics (the {\em diamond} pseudo-metric) on them generalizing the induced Riemannian metric \cite{oM-GH}. Two g.h. LLSs $L_1, L_2$ are called {\bf equivalent at} their Cauchy sets $S_1$ resp. $ S_2$ iff there are isometric  neighborhoods $U_i$ of $S_i$ in $L_i$, and we denote the class a {\bf Cauchy germ}. A {\bf (Einstein) Cauchy development} of a Cauchy germ $G$ is a Ricci-flat (in the $L^1$-optimal transport sense, w.r.t. Lorentzian Hausdorff measure) LLS in the equivalence class $G$. The pseudometric on a Cauchy set depends on its germ only, but it is unclear how to generalize the {\em second} fundamental form. 

\newpage 

Among the concepts from classical relativity whose transfer to LLSs is challenging are the following: 
\begin{enumerate}
\item There is no compact exhaustion natural in the category of pointed globally hyperbolic (g.h.) Lorentzian manifolds, not even in each Minkowski space as a single object with its morphisms: no compact subset containing another point is invariant under all boosts around the puncture;
\item Each (regular g.h.) LLS $X$ is extendible as a (regular g.h.) LLS (choose an arclength-parametrized maximizer $c:  [0; D] =: I \rightarrow X$, let $A := \{ x \in \R^{1,1} \vert x_1 \geq 0, x_1 \leq 2 x_0 \leq 2D- x_1 \} $, a triangle with one timelike and two spacelike edges, and glue $X$ to $A$ along $c(I)$ and the timelike edge $A \cap x_1^{-1} (\{ 0 \})$;
\item For local uniqueness of vacuum Cauchy developments, it would be highly desirable to have "coordinates" (i.e., a canonical representative of an isometry class of LLSs), like harmonic ones in the spacetime case, that do not depend on the choice of a Cauchy time function (they could, however, be allowed to depend on an initial Cauchy surface $S$). How to define them?
\end{enumerate}

Surprisingly, in this article, to each of the above issues we find a pleasant answer, which each time involves {\em geodesic well-posedness} (cf.Sec.1), i.e., unique extendibility of geodesics defined on compact proper intervals and $C^0$ dependence on initial values. Arguably, this assumption is natural from the viewpoint of physics, encoding that motion of particles is predictable. In Section 0, we suggest a way to "naturalize" LLSs. Part of the data of an LLS \cite{KS} is a metric, which is not natural in e.g. the category $S_n$ of $n$-dimensional spacetimes with time-oriented isometries as morphisms (cf. \cite{oM-3var}), so there is no functor from $S_n$ to any category of LLSs in which the morphisms preserve the metric: There is no natural identification of (g.h.) spacetimes with (g.h.) LLSs. However, what is needed in the applications in the foundational paper \cite{KS} is only a local bi-Lipschitz equivalence class of a metric (the topology can be defined by other means). In Sec. 0 we perform a minimally invasive surgery on the definition of LLSs, by defining {\em Lorentzian space} (LS) to allow not only for the obvious forgetful functor from the category of LLSs but also for a functor from a category of spacetimes, and address each of the three issues above in a particular section, obtaining:

\begin{enumerate}
\item There is a natural compact exhaustion in the category of {\em geodesically well-posed (g.w.)} g.h. LSs that are {\em doubly} pointed (i.e., with the additional datum of two points $p,q \in L$ with $p \ll q$), and, as an implication, a doubly-pointed Lorentzian Gromov-Hausdorff metric such that the maximal Cauchy development for spacetimes is a continuous map. By inducing natural Riemannian metrics on doubly pointed spacetimes, this could help extending the finiteness results in \cite{oM-Lorfin} from Cauchy slabs to general spacetimes and possibly LSs.
\item There is a (non-spacetime) g.h. LS properly extending Kruskal spacetime, which, {\em if geodesically well-posed}, can be proven to be maximal.
\item We will construct "Hades coordinates"\footnote{name due to the fact that we use "shadows" of points to "coordinatize" extensions of given Cauchy data}: canonical representatives of isomorphism classes of {\em geodesically well-posed} Cauchy developments of a noncompact connected Cauchy set.
\end{enumerate}

\newpage 

{\bf 0. Preliminaries on Lorentzian (length) spaces} 
 
 \bigskip
 
  A {\bf Lorentzian space} is a $4$-tuple $(X, \leq , \tau, d)$ where $X $ is a nonempty set and
 	$\leq$ is a binary relation on $X$; $\tau: X \times X \rightarrow [0; \infty ]$; $d: 2^X  \rightarrow [0; \infty]^{X \times X}$;
  satisfying the following properties P.1 to P.6: 
 
 \begin{enumerate}[{P}.1]
 	\item $\leq $ is an order relation, we define $y \in J^+(x) : \Leftrightarrow x \in J^- (y) : \Leftrightarrow y \geq x : \Leftrightarrow  x \leq y \ \forall x,y \in X$.
 	\item $\forall x,y,z \in X: x \leq y \leq z \Rightarrow \tau (x,z) \geq \tau (x,y) + \tau (y,z)$.
 	\item For a real interval $I = [a; b]$, call $c: I \rightarrow X$ a {\bf future timelike} resp. {\bf causal curve} iff there is a neighborhood $U$ of $c([a;b])$ s.t. $c$ is $[d_U]$-Lipschitz and for all $s,t \in I$ we have $s < t \Rightarrow c(s) \ll c(t) $ resp. $s < t \Rightarrow c(s) \leq c(t) $, and let $\Omega(X,x,y) $ be the set of future causal curves in $X$ from $x$ to $y$. With $y \in I^+(x) : \Leftrightarrow x \in I^-(y) : \Leftrightarrow x \ll y :\Leftrightarrow  y \gg x : \Leftrightarrow \tau (x,y) >0 $, we require $\ll \subset \leq$ and $ x \leq y  \Leftrightarrow \Omega (X,x,y ) \neq \emptyset $.
 	\item  For $A\subset X$ we define the {\bf joint future resp. past $J^{\pm \cap} (A) := \bigcap_{a \in A} J^\pm (a)$ of $A$}. We define $\liminf \ ^\pm , \limsup \ ^\pm: X^\N \rightarrow P(X)$ and $L: X^\N\rightarrow P(X)$ by, for all $a: \N\rightarrow X$:
 	
 	\bean
 	\liminf \ ^\pm (a) &:= & \bigcup_{n \in \N} J^{\pm \bigcap} (\{ a(m) \vert m \geq n \}) ,\ 
 	\limsup \ ^\pm (a) := \bigcup \{ J^{\pm \bigcap} (a \ci j (\N)) \vert j: \N \nearrow \N \},\\
 	L(a) &:=& \{ v \in X \vert I^\pm(\liminf  \ ^\pm(a)) = I^\pm(\limsup \ ^\pm (a)) = I^\pm(v) \}.
 	\eean

 	\V We define \footnote{If applied to causally simple spacetimes, the topology ${\rm tp}$ is the Alexandrov interval topology. Conversely, $(X, \leq, tp)$ is always causally simple, in particular locally causally convex, for details see \cite{oM-fcc}, Sec. 8.} a topology ${\rm tp}$ by $S \subset X$ tp-closed iff $L(a) \subset S \forall a: \N\to S$. We require $\tau$ to be lower ${\rm tp}$-semicontinuous, and tp to be metrizable and sigma-compact. For each $3$ compact sets $C,K,L$ with $L \subset C,K$ we require $d(C)|_{L \times L}$ to be a metric bi-Lipschitz equivariant to $d (K)|_{L \times L}$. In particular, the map $d $ defines a local Lipschitz class $[d]$.
  	\item $\tau $ is {\bf intrinsic}, i.e., $\forall x,y \in X: x \leq y \Rightarrow \tau (x,y) = \sup \{ \ell^\tau (c) \vert c: x \leadsto y {\rm \ causal}\}$. 
 	\item $U \subset X$ is {\bf localizing} iff $\exists u \in C^0(U \times U , [0; \infty))$ with $u \leq \tau $ s.t. for some (hence any) $d \in [d]$:
 	
 	(i) There is $C>0$ s.t. each causal curve contained in $U$ has $d$-length $<C$, 
 	
 	(ii) $\mathcal{U} := (U,   \leq\vert_{U \times U}, \ll\vert_{U \times U}, u, [d\vert_{U \times U}])$ satisfies Items P1,P2,P3,P4,P5.
 	
 	(iii)  $\mathcal{U}$ is {\bf geodesic}: $\forall x,y \in U : x <_U y \Rightarrow \exists k \in \Omega(U,x,y): \ell^\tau (k) = u(x,y)$. 
 	
 	
 	Each $x \in X$ has a causally convex localizing neighborhood, and there is a connected and dense subset $D $ of $X$ s.t. each $x \in D$ has a localizing neighborhood $U$ with $I_U^\pm (y) \neq \emptyset \forall y \in U$.
 
 	\end{enumerate}
 
If we required the relation $\ll$ to be full, we could have just used the Alexandrov topology instead of $\tau_+$, and in this way we would not have to restrict to causally simple spacetimes but could admit causally plain strongly causal $C^0$ spacetimes as in \cite{KS}. However, in view of the gluing construction in Sec. 2 we prefer not to exclude e.g. black holes with future boundary, diamonds and Cauchy slabs (see below), all violating the fullness assumption $I_U^\pm (y) \neq \emptyset \forall y \in U$ of the original definition, and thus use the more complicated topology $\tau_+$, which requires causal simplicity.

The notion of Lorentzian space is a small variation on the notion of weak LLSs in \cite{oM-3var}. Each causally simple spacetime is a Lorentzian space if equipped with local Noldus$^p$ metrics, see below. As mentioned in the last footnote, by choice of the topology tp, the space $(X, \leq, {\rm tp})$ is causally simple, thus the requirement of local causal closedness as in the definition of LLS would be redundant here. 
An LS $X$ is called {\bf globally hyperbolic}, short {\bf g.h.}, iff for each compact $C \subset X$, each causal curve $c: [0; 1) \rightarrow C $ has a continuous extension to $1$, and for each $p,q \in X$ we get $J(p,q) $ compact. This definition implies by \cite{BGH}, Th. 2.13, 2.27, that a g.h. LS is non-totally imprisoning, and then $ \tau $ is finite and continuous, and in   P.5 we can take $u = \tau\vert_{U \times U}$. An LS in $X,Y$ in which each two $x,y \in X$ with $x \ll y$ can be joined by a {\em timelike} curve is called {\bf time-connected}\footnote{Whereas time-connectedness is desirable from the point of view of classical partical trajectories and also part of the original definition of LLSs, the example in \cite{GHS} together with the definition of the Gromov-Hausdorff metric in \cite{oM-GH} shows that the space of time-connected g.h. LSs is not closed in the space of all g.h. LSs.}. An LS in which each maximizer containing a null segment is null is called {\bf regular}. By \cite{BKR}, Lemma 3.6, an LLS is regular if and only if each point has a localizable regular neighborhood. 
An LS each of whose points has a localizing neighborhood basis is called {\bf strongly localizable}. The subset $D$ as in P.5 is always open by definition, and we easily see that if $X$ is regular, each two points $p,q \in D$ can be joined by a piecewise timelike curve.

For an LS $X \ni p$ let $\de (p) := \max \{ \sup \{ \tau (p, x) \vert x \in X    \},  \sup \{ \tau (x,p) \vert x \in X   \} \}$. An LS is called {\bf temporally compact} if $J^\pm (x) $ are compact for each $x \in X$. A {\bf slab} $X$ is a g.h. connected temporally compact LS s.t. $\de^{-1} ([0; \e])$ is compact for some $\e >0$. Examples of slabs are Cauchy slabs $X$ of positive Lorentzian diameter $\inf_X \de$ and causal diamonds of temporally related points. $X$ is called {\bf strong} iff for each slab $U \subset X$, each Cauchy sequence w.r.t. $d |_{ U \times  U} $ has a limit in $\cl (U)$.

 Let $U \subset X$ be open, $p, q \in (0; \infty )$, let $\sigma (x,y) := \frac{1}{2} (\tau (x,y) - \tau (y,x))$, the antisymmetrization of $\tau$. Then $\check d^U_{p,q}: (x,y) \mapsto \sup \{ (|\sgn(\s(x,z)) \cdot \vert \sigma (x,z)\vert^p - \sgn (\sigma (y,z) ) \cdot \vert \sigma (y,z)\vert^p |)^{q} : z \in U \}$ is the {\bf $(p,q)$-Noldus metric of $U$} and $d^U_{p, q} \geq \check{d}^U_{p,q}$ the intrinsification of $\check d^U_{p, q}$, called the {\bf intrinsic $(p, q)$-Noldus metric of $U$}. For $X$ a g.h. spacetime, the $\check{d}_{p,q}^X$ are in general locally unbounded w.r.t. the manifold topology, but for any Cauchy slab and for all $p \geq 1$ the intrinsic $(p, q)$-Noldus metric is a true (finite) metric if and only if $q = 2/p$, and then its bi-Lipschitz equivariance class is the one of the Cauchy slab, cf. Th. 14 in \cite{oM-GH} \footnote{as well as the intrinsification $d_{2,X}$ of $\check{d}_{2,X}: (x,y) \mapsto || \tau (x, \cdot ) - \tau (y, \cdot) ||_{L^2 (X)} $ cf. Th. 7 in \cite{oM-3var}, which, however, transfers to LLS only in the presence of a natural measure, like a Lorentzian Hausdorff measure as in \cite{MS}.}. A Lorentzian space $X$ is called {\bf supremal} iff for each slab $U$ in $X$ there is a representative $d$ of $[d]$ and $c,C,p,P,q,Q >0$ with $cd^U_{p,q} \leq d \leq C d^U_{P,Q}$ (*) and {\bf infinitesimally Hilbertean (IH)}\footnote{Interesting examples of supremal non-IH LSs are defined in \cite{jR}.} iff we can choose $p = P , Q=q=2/p$ in (*).

 We define categories {\bf LLS} of causally simple Lorentzian length spaces, {\bf LSP} of Lorentzian spaces, the morphisms $f: X \rightarrow Y$ being in each case those pulling back every relation or map involved, e.g. $f^* \tau_Y = \tau_X$, $x_1 \leq x_2 \Leftrightarrow f(x_1 ) \leq f(x_2)$. We consider the category {\bf SCS} of causally simple Lorentzian spacetimes and causally convex time-oriented isometric embeddings.
 {\bf SLS} is defined as the subcategory of {\bf LSP} of the strongly localizing, regular, strong supremal objects. Let $F$ map an LLS with a metric $d$ to the LS with the equivalence class $[d]$ of $d$ and let $T$ map a strongly causal spacetime $(M,g)$ to the object $(M,\tau_g)$. Recall that there is no functor from {\bf SCS} to {\bf LLS}. The assignment $A$ of an LS to an element of {\bf SCS} is as follows: As a local Lipschitz structure take the manifold Lipschitz structure (even a global Lipschitz class). It takes values in {\bf SLS}: For the strongness of the image see \cite{oM-3var}, Th.6. 
 Thus we get functors $T: {\rm \bf SCS} \rightarrow {\rm \bf LSP}$, $A: {\rm \bf SCS} \rightarrow {\rm \bf SLS}$, $F: {\rm \bf LLS} \rightarrow {\rm \bf LSP}$, which satisfy $A=T$: Standard ODE estimates in local coordinates show that the local Lipschitz class of $d_{2,1}$ is the restriction of the global, manifold, Lipschitz class. The functor $F$ leaves $X, \leq, \tau$ unchanged, forgets $\ll$ and replaces the metric with the Lipschitz class of its restrictions to compacta. This suggests replacing the category {\bf LLS} with the category {\bf LSP}.

 \bigskip
 
\V Going through the proofs in \cite{KS} and \cite{BGH} and defining $\partial^\pm X:= \{ x \in X | I^\pm (x) = \emptyset \}$ we get

\begin{Theorem}
	All assertions on g.h. LLS in \cite{KS} and \cite{BGH} are true for g.h. LS. Even more, for $X$ g.h. LS we get $J(A,S) $, $J(A,B)$ compact for $A,B \subset X \setminus \partial^+ X \setminus \partial^- X$ compact, $S$ a Cauchy set of $X$. 
\end{Theorem}

 \V {\bf Proof.} The nontrivial part, showing properness of the metric in \cite{BGH}, works by: 
 
 \begin{Theorem}
 	\label{ProperMetricInLocalLipschitzClass}
 	Each local Lipschitz class of an intrinsic metric $d$ on a sigma-compact metric space $M$ contains an intrinsic proper metric.
 \end{Theorem}

\V {\bf Proof of Th. \ref{ProperMetricInLocalLipschitzClass}.} A compact exhaustion $\{ C_n | n \in \N\}$ is called {\bf proper} iff $C_n \subset \Int C_{n+1} \ \forall n \in \N$. Sigma-compact metric spaces admit proper compact exhaustions. Let $n \mapsto C_n$ of $M$ be a proper compact exhaustion of $M$. For each $n \in \N$ we define $ \e_n:= d(C_{n-1}, \cl (M \setminus C_n))$. Let $f \in C^0( M , (0; \infty))$ with $f|_{C_n} \geq (\inf \{ \e_k | k \in \N_n \} )^{-1}$. Replace the length functional $\ell= \ell^d$ with another length functional $\tilde{\ell} (c)= \sup \{ \sum f(c(t_i)) \cdot d(c(t_i), c(t_{i+1})) \vert t \in P(I) \} $ for each $c: I \rightarrow M$, where $P(I)$ is the set of partitions of $I$. The metric $D$ generated by $\tilde{\ell}$ is intrinsic by definition and satisfies $D(C_1, M \setminus C_n) \geq n-1$: Let a locally Lipschitz curve $c : I \rightarrow M , c: C_1 \leadsto C_n $ be given, let $c(t_k)$ be the first point in $C_k$ for each $k$, then $\ell (c) \geq \sum_{k=1}^{n-1} D(c(t_k), c(t_{k+1}) \geq \e_k^{-1} d(c(t_k), c(t_{k+1})) \geq n-1 $. In particular, $D$ is proper. \hfill \qed
 
\bigskip

Consequently, we can choose a proper representative in the given local Lipschitz class and apply the theorems proven in \cite{BGH}. For the last assertion, mind that in the proof of the limit curve theorem in \cite{BGH} we just need $c_n(t) \rightarrow_{n \rightarrow \infty} x_0$ which we have for $c_n(t) $ contained in a compact set. \hfill \qed

\newpage 
 
{\bf 1. Natural compact exhaustion in LS, pointed Gromov-Hausdorff metrics}

\bigskip

Let us define a synthetic version of well-posedness of the geodesic equations. For a $C^2$ Lorentzian metric, the Christoffel symbols are $C^1$, therefore the geodesic equations are locally Lipschitz and thus {\em well-posed}, i.e., in any spacetime $(M,g)$, given $v \in TM$ we have local existence and uniqueness of a geodesic $c_v$ with $c_v'(0)= v$, and $c_v$ depends continuously on $v$. In the absence of a tangent space for a general LS, we now rephrase well-posedness in the new context: An LS $L$ is called {\bf non-branching} iff there is no branching timelike geodesic in $L$, i.e., any two local maximizers agreeing on a nonempty open interval agree up to parametrizations and restrictions. Each LS with timelike sectional curvature locally bounded below is non-branching, \cite{BS}. Moreover, Prop. 4.6 in \cite{GKS} implies that in each g.h. LS a geodesic $[0; T) \rightarrow X$ is extendible to a geodesic on $[0;T]$ if and only if it is extendible as a causal curve. An LS $X$ is called {\bf non-stopping}\footnote{In \cite{BS}, Def. 4.2, this is called the geodesic prolongation property.} iff each causal geodesic (i.e., locally maximizing) curve $c: [0;T] \rightarrow X$ with $I^+(c(T)) \neq \emptyset $ has a geodesic extension to a strictly larger interval $[0; S) $, and time-dually. From now on assume $X$ regular, implying the existence of parametrization by arc length. For all arc-length parametrized geodesics $c: [a_-;a_+] \rightarrow X$ in a non-branching and non-stopping LS $X$ define $E(c): (T^-(c) ; T^+ (c) ) \to X$ as the maximal arc-length parametrized geodesic extension of $c$ (existence implied by \cite{BGH} (2.13, 2.20, 2.22)). For $\e>0$, define 

\bean
 T_\e^- (c)  &:=& \max \{ a_- - (a_+-a_-) , T^- (c) + \e (a_- - T_\e^-(c))\} , \\ T_\e^+ (c) &:=& \min \{a_+ + (a_+-a_-) , T^+(c) - \e (T_+(c) -a_+) \}, \\ E_{ \e} (c) &:=& E (c)|_{[T_\e^- (c); T^+(c)]} .
 \eean

 In words, $E_{\e}$ triples the length of geodesics as long as their parameter stays a portion determined by  $\e$ away from the boundary of the domain of definition. Let $R_\tau$ be the reparametrization of timelike curves by $\tau$-arclength, and for $p,q \in X $ and for a causal curve $c$ in $J(p,q)$ let $R_{p,q} (c)$ be the reparametrization of $ c$ according to the time function $\tau_{pq} := \tau (p, \cdot ) - \tau (\cdot, q)$. An LS $X$ is called {\bf geodesically continuous} if and only if it is regular and for each geodesic $c: [a;b] \rightarrow \Int X$ and each $\e >0$ there is an open neighborhood $U  $ of $c |_{[a+ (b-a)/3 + \e; b -  (b-a)/3 - \e]}$ in the space $M(c(a + (b-a)/3), c(b - (b-a)/3))$ of maximizers in $J(p:= c(a + (b-a)/3), q:= c(b - (b-a)/3))$ parametrized by $\tau_{p,q}$, equipped with the metric $d_{J(p,q)}$, such that for all $k \in U$and all $\e >0$  we have ${\rm im} (E_\e (R_\tau (k))) \subset J(c(a), c(b))$ and $R_{c(a), c(b)} \ci E_\e  \ci R_\tau|_U $ is continuous\footnote{Whereas this definition may appear to be over-complicated, it turns out to be essentially the only reasonable definition of continuity of geodesic extension in the absence of a global metric.}.  
$X$ is called {\bf geodesically well-posed (g.w.)} iff it is non-branching, non-stopping and geodesically continuous. In g.w. spaces a causal geodesic is extendable as a causal geodesic if and only if it is extendable as a causal curve. $X:= \{ x \in \R^{1,n} | x_1 \geq 0\}$ is non-branching but stopping: there are timelike geodesics stopping on the timelike boundary. $X$ has vanishing timelike sectional curvature. However, in $X$ we have {\em horismos with maximal points}, i.e. there are $x \in X$ and $y \in E^+(x) := J^+ (x) \setminus I^+(x)  $ with $ J^+(y) \setminus \{ y\} \subset I^+(x) $. Conversely, we can construct an LS $Y$ conformally related to $\R ^{1,1}$ with a non-continuous conformal factor $\omega$ (s.t. $Y$ has no horismos with maximal points) that is non-branching but stopping: Define $\omega (x) = 1+ |x_1|^{-1/2} \cdot f(\theta (x) )$ for all $x \in Y \setminus \{ 0\}$ with $\theta (x)$ being the angle of $x$ to $\R^+ \cdot \partial_1$ and $f $ differentiable nonnegative with $f ([\pi/4; 3 \pi/4])= f([5 \pi/4; 7 \pi/4]) = 0$ and $f ([7\pi/8; 9 \pi/8])= f([15 \pi/8; 2 \pi) = f([0; 1 \pi/8] = 1$. Then $Y$, which is even a Gromov-Hausdorff limit (see \cite{oM-GH}) of spacetimes all conformal to $\R^{1,1}$, has no geodesic through $0$ and has unbounded curvature at $0$.

\begin{Theorem}
	Globally hyperbolic Lorentzian spacetimes are geodesically well-posed.
\end{Theorem}

\V{\bf Proof.} The only non-standard part is proving geodesic continuity, as the neighborhood is taken in the $C^0$ topology, not in the $C^1$ topology. But as the neighborhood is in the space of geodesics, the topologies are equivalent (seen in a geodesically convex neighborhood of the initial point $c(a)$). \hfill \qed

\bigskip

Let $X$ be an LS and let $\e >0$.

For a maximal curve $c$ in $X$, let $p^\pm_\e (c) := E_{\e}(c)(T^\pm_\e (c))$. Define maps $Z_{\e}^\pm, A_{\e}: 2^X \rightarrow 2^X$ as follows:

\[\forall \e>0, K \in 2^X: \de(K, \e):= \{ (p,q) \in K^2 | \tau (p,q) \geq \e\}  \subset X \times X, \]

\[z_{\e}^\pm (p,q) :=  \bigcup \{ p_\e^\pm (c) | c: p \leadsto q {\rm \ maximal} \} , \ Z^\pm_\e (K) := \bigcup_{(a,b) \in \de(K, \e/3)} z^\pm_\e (a,b),  \ A_{\e}(K) := J(Z_{\e}^- (K) , Z_{\e}^+ (K)). \]

for $K \in 2^X$. For $(p,q) \in X^2$ with $p \ll q$, let $\e (p,q):=  \frac{1}{3}\min\{ \tau (p,q), 1\} >0$ and $C_n (p,q):= (\bigcirc_{k=1}^n A_{\e (p,q)/k}) (J(p,q))$, where, for an $n$-tuple of maps $f_i: T \rightarrow T$, let $\bigcirc_{i=1}^n f_i$ be the iterated application $f_n \ci ... \ci f_1$ of maps, defined inductively. Let $C(X)$ be the set of compact subsets of $X$. 


\begin{Theorem}
\begin{enumerate}
	\item $\de(K, \e) \in C(X \times X)$ for $K \in C(X)$, and $Z_\e^\pm, A_{\e}: C(X) \rightarrow C(X)$.
	\item $\{ C_n (p,q)| n \in \N \} $ is a proper compact causally convex exhaustion for each connected g.w.g.h. LS $X \ni p,q$. It is natural in the category of doubly-pointed g.w.g.h. Lorentzian spaces.
\end{enumerate}	
\end{Theorem}

\V{\bf Proof.} For Item 1: Continuity of $\tau$ implies the claim on $\de ( \cdot , \e) $. We show the assertion on $Z_\e^\pm$, the one on $A_\e$ follows. So we want to show compactness of $Z^\pm_\e (K)$ for $K$ compact. For all $p \in X$, the function $t := \tau ( p , \cdot )$ is a time function on $I^+(p)$. Each maximal curve from $p$ to $q \gg p$ parametrized by arc length is parametrized by $t$. The space $Y$ of such is compact by a Limit Curve Theorem due to Minguzzi (Th. 10 in \cite{eM-Results}). Thus we choose finite sets $\{ p^\pm_1,... p^\pm_N \} $ with $ \de(K, \e/3) \subset \bigcup_{j=1}^{N} I (p_j^-, p_j^+)$. Then $A:= \bigcup \{ \Omega (p,q) |  (p,q) \in \de^- (K, \e/3) \}  $ is a (closed) subset of the (compact) set $\bigcup_{l=1}^{N} \Omega (p_l^-, p_l^+)$. So $A$ is compact. For each element $c$ of $A$, we can choose its neighborhood $N_c:= J(p^-_{\e/2}, (c), p^+_{\e/2} )$. Compactness of $A$ implies that $A$ is covered by a finite number of those $N_c$. Geodesic continuity shows that the set of endpoints of the extensions in each $N_c$ is compact.

For Item 2: First we show properness. Let $x \in C_n(p,q) $, then by the recursive definition, there are two tuples $(a^\pm, b^\pm)  \in \de (C_{n-1} (p,q), \frac{\e}{3n})$ and $p^\pm \in z_{\e/n}^\pm (a^\pm, b^\pm)$ such that $x \in J(p^-,p^+)$. Thus there are future timelike curves $c: p^- \leadsto a^- \leadsto b^-$ and $ c^+: a^+ \leadsto b^+ \leadsto p^+$ with $\ell (c^\pm) > \min \{ \tau (a^+, b^+, a^-, b^- ) \geq \e(p,q)/n$. The triangle inequality then implies that $(p^-, p^+) \in \de (C_n , \e/n) \subset \de (C_n, \e/(n+1) ) $. So properness of $C(p,q)$ follows from 

$x \in J(p^-,p^+) \in I(Z^-_{\e(p,q)/(n+1)} (C_n(p,q)), Z^+_{\e(p,q)/(n+1)} (C_n(p,q))) \subset \Int (C_{n+1} (p,q))$. 


For the exhaustion property, let $x \in X$, we want to show that there is $n \in \N$ with $x \in C_n$. Let $c: [0;1] \rightarrow X$ be a continuous (i.g. noncausal) curve from $y \in I(p,q)$ to $x$ and let $T:= \sup \{ x \in [0;1] | \exists m \in \N: c(s) \in C_m (p,q)\}$ ($>0$, as $C_1 (p,q)\supset I(p,q)$ is a neighborhood of $y$). Let $a,b \in X$ with $c(T) \in I(a,b) $. Then there is $s <T$ with $c(s) \in I(a,b)$. There is $m \in \N$ with $c(s) \in \Int (C_m (p,q))$. Let $k^+ := E(\gamma) |_{[0; T_+(\gamma)]}$ be the maximal future extension of a maximizer $\gamma $ from $c(s) = \gamma(0)$ to $b$ as a geodesic parametrized by Lorentzian arclength, analogously we define $k^-, T_-(\gamma)$. Then $k^\pm$ intersects $\partial^\pm C_m (p,q)$ at $z^\pm $. Defining $T_1, T_2>0$ by $ k^+(T_1 )= z^+$, $k^+(T_2) = b$, we get $c(s), z^+ \in C_N(p,q)$ for all $N \geq m$. Now choose $N^+ \geq m$ with $ N^+ \epsilon (p,q) \geq \frac{T_2 - T_1}{T_1}, \frac{T_2 - T_1}{T_+(\gamma)- T_2}$. Then $b \in E_{\epsilon (p,q)/N^+} ((T^-_{\epsilon (p,q)/N^+}, T^+_{\epsilon (p,q)/N^+}))$. By time-dual arguments we get $ N^- \in \N$ with $a \in E_{\epsilon (p,q)/N^-} ((T^-_{\epsilon (p,q)/N^-}, T^+_{\epsilon (p,q)/N^-}))$. If $T<1$, for $N := \max \{N^-, N^+ \}$ there is some $\rho>0$ with $ c(T+ \rho) \in I(a,b) \subset C_N(p,q) $, contradiction, thus $T=1$, $x \in C_N(p,q)$. \hfill \qed

\bigskip

Define the distance of LSs $(X_1, \tau_1)$ and $(X_2, \tau_2)$ w.r.t. subsets $A_1 , A_2 \subset X_1$ and $B_1, B_2 \subset X_2$ by 

\vspace{-0.4cm}

\[ {\rm Corr}_{A, B} (X,Y) := \{ \rho \in {\rm Corr} (X,Y) | \rho \cap (A_i \times X_2) \subset A_i \times B_i \supset \rho \cap (X_1 \times B_i)  \ \forall  i \in \{ 1,2 \} \} ,\]   

\vspace{-0.3cm}

\[ d^- ((X_1, A), (X_2,  B)):= \inf \{ {\rm dis} (\rho  | \rho \in {\rm Corr}_{A,B} (X_1,X_2) \}  \in [0; \infty] \]

where ${\rm dis}$ denotes the distortion of $\rho$. Note that $d^- ((X_1, A),(X_2, B)) < \infty$ for $A,B$ compact. For $A_i = \{ (p_i, q_i) \in J_i \}$ we define a Gromov-Hausdorff metric for doubly-pointed LSs by 

\vspace{-0.2cm}

\[ d^- ((X_1, p_1, q_1), (X_2, p_2,q_2)) :=  \sum_{i=0}^\infty 2^{-i} \arctan (d^-(  K_i ((p_1, q_1) , \{ p_1, q_1\}), ( K_i (p_2, q_2) , \{ p_2, q_2 \} ))) .\]

By the above natural compact exhaustion we can transfer everything done for (spacetime) Cauchy slabs in \cite{oM-GH} to (i.g. spatially and temporally noncompact) g.w. g.h. regular LSs. Let $d_{N, K_m} $ be the Noldus metric w.r.t. $K_m$, $d_{N, K_m} (x,y) := \sup \{ | \tau^2(x,z) - \tau^2 (y,z) |:  z \in K_m\}$, then on a doubly pointed LS $(X,p,q)$ define a metric $d_{p,q}$ by \(d_{p,q} (x,y) := \sum_{m=0}^\infty 2^{-m} \arctan (d_{N,K_m} (x,y)) \ \forall x ,y \in X\).

\begin{Theorem}
$d_{p,q}$ is natural in the category of g.w. g.h. regular doubly punctured LSs and generates the Alexandrov interval topology on $X$. If $X$ is a spacetime and $S$ is a $C^1$ hypersurface, $d_{p,q}|_{K \times K}$ is bi-Lipschitz equivalent to the induced Riemannian metric on $K$ for each compact $K \subset S$.
\end{Theorem}

\V{\bf Proof.} In \cite{oM-GH}[Th.14] it is shown that the Noldus$^2$ metric w.r.t. a compact set is locally Lipschitz continuous. The converse estimate is directly seen in normal coordinates. \hfill \qed 
 
\newpage
 
{\bf 2. A non-spacetime maximal g.h. LS extending Kruskal spacetime}
 
 	\bigskip
 	
For each LS $Y$ let $\partial^\pm Y$ be its future resp. past causal boundary, which equals the set $TIP(Y) $ resp. $TIF(Y)$ of indecomposable past resp. future subsets of $Y$ (see \cite{oM-fcc}, where this set is given a metrizable topology ${\rm tp}$) and assume that $\tau$ has a continuous extension $\hat{\tau}: \hat{Y}^\pm \times \hat{Y}^\pm \rightarrow \R$ where $\hat{Y}^\pm := Y \cup \partial^\pm Y$. (At the end of this section, we will see that this is the case for $Y^\pm$ being the interior Schwarzschild solution resp. its time-dual.) For $U \subset X$ compact, this agrees with the definition $\partial^\pm U:= \{ x \in \partial U | I^+(x) \cap U = \emptyset\}$. Let $C(T)$ be the set of compact subsets of a topological space $T$.

For $ x \in \hat{Y}^+$ we define $G^-(x) := I^-(x) $ for $x \in Y$ and $G^-(x) = x $ for $x \in \partial^+ Y$. We define the LS structure on $X_+:= \hat{Y}^+$ by

	\begin{enumerate}
	\item extending $\leq_Y$ to $\leq_X \subset X \times X $ by $u \leq_X v : \Leftrightarrow G^-(u) \subset G^-(v)$ for all $u,v \in \hat{Y}$, 
	\item extending $\tau_Y$ to a function $\tau: X \times X \rightarrow [0; \infty) $ by $\tau(u,v) := \inf_{w \in G^-(u)} \sup_{z \in G^-(v)} \tau (w,z)$,
	\item defining $d_H$ as unique continuous extension of $d_{H \cap Y}$ for all $H \subset X$ w.r.t. the topology tp.
\end{enumerate}

We perform all these operations time-dually to induce the structures on $\hat{Y}^-$.

 		Let $(X_{\pm }, \leq_{\pm }, \ll_{\pm }, \tau_{\pm }, [d]_{\pm })$ two g.h. LSs, let $A^-  \subset X_- $ be future and $A^+ \subset X_+ $ be past (this is satisfied in particular if $A^\pm  := \partial^\mp X_{\pm } \subset X_\pm$ is a Cauchy set of $ X_\pm$, in which case $\tau|_{A^\pm \times A^\pm} = 0$). Assume $ \forall C \in C(X_+): J^-(C)  \in C(A_+)$ or $\forall C \in C(X_-): J^+(C)  \in C(A_-)$, and that there is a homeomorphism $F: A^- \rightarrow A^+$ with $F^* \tau_+ = \tau_-|_{A^- \times A^-}$. Then we call the tuple $(X_\pm, A^\pm, F)$ a {\bf gluing set} and define the {\bf gluing $X(X_{-}, F, X_+)$ of $X_-$ to $X_+$ through $F$} by

 	\begin{enumerate}
 		\item the set $X = (X_{-} \cup X_+)/ (A^- \ni x \sim F(x) \in A^+)$ (yielding injective maps $h^\pm : X_{\pm} \rightarrow X$), let $A_0:= h^\pm (A^\pm ) = [A^-] = [A^+]$,
 		\item $\leq_X \subset X \times X $ by $\leq_X := \leq_{X_{-}} \cup \leq_{X_+} \cup \{ (x,y) \in X_{-} \times X_+: \exists g \in A_0 : x \leq_{-} g \leq_+ y \}$,  
 	
 		\item $\tau_X: X \times X \rightarrow [0; \infty) $ by (for $p,q \in X$) defining $D(p,q) := J(p,q) \cap A_0$, $\tau_X |_{X_{\pm } \times X_{\pm }} = \tau_\pm$ and $ \tau_X (p,q) := \sup \{ \tau_{-} (p, d) + \tau_+ (F(d), q ) | d \in D (p,q) \}$ for $D(p,q) \neq \emptyset$ and $0$ otherwise,
 		\item $d_X: 2^X \rightarrow [0; \infty]^{X \times X} $ by $d_X (H) (p^-,q^+) := \inf \{ d (H \cap X_{-} ) (p,r) + d(H \cap X_+) (F(r),q)  | r \in A^- \}$, and s.t. $d_X (H)$ extends $d^- (H \cap X_-) $ and $d^+ (H \cap X_+)$. 
 		\end{enumerate}
 	
Note that the so defined Lorentzian distance functions are indeed intrinsifiable.\footnote{One could even admit $A^\pm$ to be a connected component of the space $\partial^\mp_f X^\pm$ of ideal boundary points at finite distance. Then closedness of $A^\pm$ in $\partial^\mp_f X_{\pm 1}$ ensures Hausdorffness of $X$, the openness ensures continuity of the Lorentzian distance function. There is another obvious variant where, if the local metrics for $X_\pm$ are intrinsic, the new metrics are intrinsic, by combining two length functions}

 	\begin{Theorem}
 	Let $(X_{\pm }, A^\pm, F)$ be a gluing set. Then $(X, \leq_X, \tau_X, [d_X])$ as above is a g.h. LS, and there are canonical embeddings of $X_-, X_+$ into $X$.
 	\end{Theorem}

 	\V {\bf Proof.} As noted above, in g.h. LSs, we can employ two modifications:
 	
 	\begin{enumerate}
 		\item In (P.3), we impose $\tau $ continuous.
 		\item In (P.5), we impose $U$ causally convex and $u = \tau|_{U \times U}$.
 	\end{enumerate}
 	
 Moreover, $X$ g.h. means that all causal diamonds $J(p,q) $ are compact and that $X$ is non-totally imprisoning, i.e. for each $K \subset X $ compact there is $D>0$ such that for each $c: \R \rightarrow K$ causal we have $\ell_d(c) <D$.   
 
 The relation $\leq_X$ is obviously reflexive, and an easy distinction of the cases $x ,y,z \in X_{-}$, $x ,y \in X_{-} \land z \in X_+$, $x \in X_{-} \land y,z \in X_+$ and $x,y,z \in X_+$ shows that it is also transitive. The statement ${\ll_X} \subset  {\leq_X}$ is true almost by definition.
 
 For the conditional triangle inequality, use again the above case distinction.  
 
 $X$ is sigma-compact as $X_{\pm }$ is: Let $n \mapsto C_n^{\pm}$ be an increasing compact exhaustion for $X_{\pm }$ then $n \mapsto K_n:= J(C_n^- , C_n^+) $ is an increasing causally convex compact exhaustion of $X$, see below.
 
 From now on assume, w.l.o.g., $ \forall C \in C(X_+): J^-(C)  \in C(A_+)$.
 
 For the continuity of $\tau_X$, let w.l.o.g. $x \in X_{-}$, $y \in X_{+}$, $x \ll y$. Then $H := J^-(y) \cap A_0 \cap J^+(x)$ is compact and nonempty. If $z_n \rightarrow_{n \rightarrow \infty} y$ then $G_n:= J^-(z_n) \cap A_0 \cap J^+(x)$ are eventually nonempty compacta Hausdorff-converging to $H$ with $n$, and 
 
 $\sup \{ \tau (x,g) + \tau (g, z_n) | g \in G_n \} \rightarrow_{n \rightarrow \infty} \sup \{ \tau (x,g) + \tau (g, y) | g \in  H \} .$
 
 For the property P.4 (intrinsicness), let wlog $x \in X_{-}$, $y \in X_+$, let $\e>0$, then there is $g \in A_0$ with $\tau (x,y) > \tau (x,g) + \tau (g, y) - \e/3$. By instrinsicness of $X_{\pm }$ there are $c_{-1}: x \leadsto g$, $c_1: G \leadsto y$ with $\ell (c_{-1}) > \tau (x,g) - \e/3$, $\ell (c_1) > \tau (g,y) - \e/3 > \tau (x,y) - \e$.

 Closedness of $\leq$ holds as $X_{-}$ and $X_+$ are closed in $X$. The first condition in our definition of 'g.h.' holds as both $J^-(q) \cap A_0$ is compact, and finally, causal diamonds $J(p,q) $ for $p \in X_-$ and $q \in X_+$ are compact: Let w.l.o.g. $\forall C \in C(X_-): J^+(C)  \in C(X_-)$, then $J^+_{X_-} (p) \cap A_-$ is compact and $J_X (p,q) \cap X_+ $ is contained in the emerald $ J_{X_+} (F(J_{X_-}^+ (p) ), q )$ and thus compact.  \hfill \qed

 	\bigskip

 	This gluing construction is closely related to the quotient time separation introduced in \cite{BR}.
 	
 	\medskip

 	Now let $K_m$ be Kruskal spacetime of mass $m>0$. Let $D^\pm :=  \partial^\pm_s K_m $ the singular part of future (resp. past) boundary defined as the set pasts (resp. futures) of $C^0$-inextendible future (resp. past) timelike curves $c: \R \rightarrow K_m$ with of finite length, or, equivalently, $(r \ci c)(t) \rightarrow_{t \rightarrow \infty } 0$ for the Schwarzschild radial parameter $r$. Let $U_m^+ $ be the black hole region, $ U_m^- $ the white hole region of $K_m$, and $R$ a (time-orientation reverting) isometry of $K_m$ swapping $U_m^+$ and $U_m^-$. It is well-known that one cannot extend $K_m$ as a Lorentzian spacetime, not even in the $C^0$ sense \cite{jS2}. We can extend the Lorentzian distance in $U_m^\pm$ continuously to their future resp. past causal completions $\widehat{U_m^\pm}$, with future and past boundary $D^\pm$, which is each a hypersurface (cf. \cite{oM-fcc}). We define 
 	

	$V_m^+:= X(\widehat{U}_m^+ , P(R) |_{TIP(U_m^+)} , \widehat{\tilde{U}_m^-}) $ (and a canonical isometric partial map $\iota^+ $ from $K_m $ to $V_m^+$), 
 	
 	where, for a map $F:K \rightarrow L$ between sets, $P(F): P(K) \rightarrow P(L)$ the corresponding map on the power set. If applied to $R$, this maps the TIPs to the TIFs. Analogously we define the set $V_m^-$, which is just a copy of $V_m^+$, with associated map $\iota^-$. This is the particular case of a gluing where $\tau|_{A_i \times A_i} = 0 $, as the $A_i $ are achronal.
 	
 	Finally\footnote{One could also replace the four subsequent gluings above by the two gluings $Y_1 := X(\hat{K}^+,  P(R), U_m^- )$ and $Y_2 := (\hat{Y_1}^-, P(R), U_m^+  )$}, consider the isometric gluing $L_m:= V_m^- \cup K_m \cup V_m^+/(\iota^- \cup \iota^+)$ (taking into account that $K_m \cap V_m^\pm = U_m^\pm $). We denote the resulting Lorentzian distance function by $\tau_L$.

 	An LS $L$ is called {\bf synoptic} or {\bf indecomposable} iff $J^\pm(x ) \cap J^\pm (y) \neq \emptyset \forall x,y \in L$ (indecomposability refers more to the equivalent property that in every decomposition $L= U \cup V $ in past subsets $U$, $V$, we have $U \subset V$ or $V \subset U$, see e.g. \cite{oM-fcc}).

 	For two LSs $L, L'$ with (equivalence classes of) germs of Cauchy sets $S, S'$ we define $L \leq L'$ iff there is an isometric embedding $i: L \rightarrow L'$ with $i(S) = S'$. Furthermore, $L$ is called {\bf maximal Cauchy development} or {\bf mCd} {\bf of $S$} iff $L \geq L'$ for each Cauchy development, and $L$ is called {\bf weakly mCd} iff for each Cauchy development $L'$ with $L' \geq L$ we have $L' = L$.

 In this context, it is a first step to look at examples of LSs that are weakly maximal as g.h. LSs. Such an example is provided by the $L_m$ defined above:

 	\begin{Theorem}
 		\begin{enumerate}
 			\item $(L_m, \tau_L) $ is a synoptic g.h. Lorentzian (length) space extending $K_m$.
 			\item If $(L_m, \tau_L)$ is geodesically non-branching, then $(L_m, \tau_L) $ is weakly maximal among the globally hyperbolic Lorentzian (length) spaces extending Kruskal spacetime.
 		\end{enumerate}
 		
 	\end{Theorem}

\V{\bf Proof.} Synopticity of $L_m$ follows from $I^\pm (p) \cap U_m^\pm \neq \emptyset$ for all $p \in K_m$ and from synopticity of $V^\pm$.

For weak maximality, let $L'$ be a g.h. LS extending $L_m$, then let $x \in L' \setminus L_m$.
Necessarily, $x $ is in the future or past of $L_m$ (as $J(x) \cap S \neq \emptyset$). W.l.o.g. $x \in J^+(L_m)$. As $L'$ is a g.h. LS, there is a future maximizer $c$ from some $y \in L_m$ to $x$ leaving $L_m$ at some $q \in \partial(L_m, L')$, $q = c(t) $ for  $t := \sup \{ s \in \R | c(s) \in K_m \}$. Let $k := c\vert_{[0;t)}$. $k$ is inextendible in $L_m$, otherwise the extension in $L_m$ and $c$ provide a contradiction to geodesic non-branching. This means that $I^-(k)$ is a TIP, either a TIP of $K_m$ or one of the additional TIPs induced by inextendible future timelike curves in $\tilde{U}_m^+$. In either case, $J^- (c(t))$ contains some part of future null infinity, and from every point of future null infinity there are past curves to $S$ of arbitrary length --- impossible if $S$ is a Cauchy surface of $L'$, and $J^-(x) \cap S $ is compact.  \hfill \qed
 	
\bigskip 	
 	
To show that $L_m$ is geodesically non-branching, regularity of the Lorentzian distance function at the future boundary of $K_m$ plays a key role. This is work in progress.	
 	
\newpage

 {\bf 3. Canonical representatives ("Hades coordinates") for developments}
 	
 	\bigskip
 	
 	A crucial step in proving existence of maximal Cauchy developments in \cite{CBG} was local uniqueness, via harmonic coordinates. A possible replacement for Lorentzian spaces could be Lorentzian strainers. Here is another way for developments of fixed initial value $S$: Let $C (S)$ be the set of compact subsets of a metric space $S$, the topology $\theta$ on $C(S)$ induced by Hausdorff distance. The maps $E_S: X \to C^0(S), X \ni x \mapsto \s(\cdot, x)|_S \in C^0 (S)$ ($\s$ is the antisymmetrization of $\tau$ as in Sec.1) and $ A_S:= \supp \ci E_S : p \mapsto J(p) \cap S: (X, {\rm tp} ) \rightarrow (C(X), \theta)$ are continuous (\cite{oM-fcc}, and \cite{ACS} Th.2, Lemma 3.11, Prop. 3.12 for $\cl (I) = J$). Now $A_S(p)$ is i.g. not enough to identify $ p \in J^+(S)$ but $E_S$ is:

 	\begin{Theorem}
 		\label{Catcher}
 	Let $X$ be a g.h. g.w. LS, let $S$ be a noncompact Cauchy set of $X$ without cut points\footnote{A cut point is a point $x$ in a topological space $X$ s.t. $X \setminus \{ x \}$ is connected. Here it implies $\# (\partial A_s (p))>1 $.}.
 	\begin{enumerate}
 	\item $ A_S^+:= A_S|_{J^+(S)}$ is increasing and $ \forall p,q \in J^+(S): A_S(p) = A_S(q) \Rightarrow p \in J(q) \setminus I(q)$. \\
$A_S^+$ is in general not injective, even for g.h. spacetimes. \\
  $A_S^+(p) \subset \Int A_S^+(q)$ does in general not imply $p \leq q$, even for g.h. spacetimes.
 	\item The map $E_S$ is strictly increasing and injective.
 	\end{enumerate}
 	\end{Theorem}

 	{\bf Remark.} The hypothesis of $S$ noncompact is indispensable, consider a Lorentzian cylinder $\R \times \mathbb{S}^1 $. The hypothesis of being geodesically well-posed (implying e.g. that $X$ is regular) is indispensable, too: modify the example in the second item in the first list of the introductory section (Sec. 0) by gluing two isometric triangles $D_1$ and $D_2$ to the same maximizing segment, then for the images $a$ resp. $b$ after gluing of the points $A:= (3D/4, D/2) \in D_1$ resp. $B=(3D/4, D/2) \in D_2$ have $A_S(a) = A_S(b)$ for all Cauchy sets $S$ with $ c(0) \in J^+(S) $.

 	\medskip
 	
 	\V{\bf Proof.} 	Let $(C,g)$ be an infinite one-ended cigar (a Riemannian manifold isometric to $\mathbb{S}^1 \times \R^+ $ outside of a compact set $K$) and $X:= (\R\times C, -dt^2 + g_C)$. There are $p, q \in X$ with $t(p) = t(q)$ and a $t$-level Cauchy surface $S$ such that $A_S(p) \subset \Int A_S(q)$ or vice versa. In $(C,g)$, $A_S$ is noninjective: $A_S ((t+ \pi, r , s)) = A_S(t,r+\pi,s)$ for $S$ an appropriate $t$-level Cauchy surface.

 	Assume $A_S(p) = A_S(q)$, then $J^-(p) \cap I^-(S) = J^-(q) \cap I^-(S) =:B$. Even more, $\partial B \setminus \int A_q(x) \subset (J^- (p ) \setminus I^- (p)) \cap J^- (q) \setminus I^+ (q))$. Pick $b \in \partial B \setminus S$ and a (null) maximizer $c_p: [0; D] \rightarrow X,  x \leadsto p$. Let $ c(T) =p $ and define $t:= \sup \{ s \in [0; D] | c_p (t) \in J^-(q) \}$. If $t < T$ and $c(t) \neq q$, then $c(t) \in  (J^- (p ) \setminus I^- (p)) \cap J^- (q) \setminus I^+ (q))$, and we get a contradiction to non-branching composing a piece of $c_p$ with a null maximizer from $c(t) $ to $q$. Thus $p \leq q$ or $q \leq p$, let w.l.o.g. the first be the case. If $p \in J^-(q) \setminus I^-(q)$, then each null line $x \leadsto p$ composed with each null line $ c_{pq}: p \leadsto q$ is maximal. Now choose $x,y \in \partial A_S(p)$ different and let $c_x : x \leadsto p$ and $c_y: y \leadsto p$ be null maximizers, then composing them with $c_{pq}$ contradicts non-branching again. The remaining case is $p \ll q$. Then (by push-up) $J^-(p) \cap S \subset I^-(q) \cap S \subset J^-(q) \cap S \subset J^- (p) \cap S$, thus the compact set $J^-(p) \cap S$ equals $I^-(q) \cap S$, which is open in $S$. As $S$ is connected, both equal the noncompact set $S$, contradiction. If $E_S(p) = E_S(q)$ then $A_S(p) = \supp E_S(p) = \supp E_S(q) = A_S(q)$. Thus w.l.o.g. $p \leq q$. Push-up implies $E_S(q) (x) > E_S(p) (x) \ \forall x \in I^-(p) \cap S$, unless $p = q$. \hfill \qed

 	\newpage

 	Let $U$ be a set, then $K \subset U \times U$ is called {\bf diagonal} iff $K= A \times A  $ for some $A \subset U$. For an LS $X$ extending $S$, we define $ A:=E_S(S) \subset C^0 (S)$. Then we can interpret the Lorentzian distance as a function $A \times A \to \R$. Let ${\bf DPR} (S)$ be the category of real functions on diagonal subsets of $C^0(S) \times C^0(S)$ s.t. the diamond metric (cf. \cite{oM-GH}) is the given one, with extensions (right inverses of restrictions) as morphisms. For a length space $S$ let ${\bf GLS} (S)$ be the category of g.h.g.w. LSs with Cauchy set $S$, with causally convex Lorentzian isometric embeddings fixing $S$ as morphisms.
 
 	\begin{Theorem}
 	Let $S$ be a noncompact connected length space. $\Phi_S: (X, \tau) \mapsto (S,(E_S^{-1})^* \tau)$ is an essentially injective functor ${\bf GLS} (S) \rightarrow {\bf DPR} (S)$, i.e. if there is an isomorphism $H: \Phi_S (X) \to \Phi_S (X') $ of ${\bf DPR} (S)$ then there is an isomorphism $h: X \to X'$ of ${\bf GLS} (S)$. 
 	\end{Theorem}	
 	
 {\bf Remark.} As there is no nontrivial isomorphism in ${\bf DPR} (S)$, this boils down to the fact that for $ \Phi_S (X) = \Phi_S (X')$ there is an isomorphism $H: X \to X'   $.

 \V{\bf Proof.} Let $X_1$, $X_2$ be two extensions of $S $ with $ E_S^{X_1}(X_1) = E_S^{X_2}(X_2)$. Then $ G:= A_S(X_1) = A_S(X_2) \subset C^0(S)$, and $(E_S^{X_2})^{-1} \ci E_S^{X_1}: X_1 \rightarrow X_2$ is an isometry of LSs. \hfill \qed

\bigskip

{\bf Acknowledgement:} The author acknowledges very helpful remarks of an anonymous referee on a first version. 

\bigskip

{\bf Funding and interest declaration.} No funding was received for conducting this study. \\ The author has no financial or proprietary interests in any material discussed in this article.

\bigskip

{\bf Data availability statement:} No experimental data has been produced for this article.

\end{document}